\documentclass[english]{amsart}
\usepackage[english]{babel}

\usepackage[hidelinks]{hyperref}

\usepackage{verbatim}
\usepackage{amsfonts,amssymb,mathtools}
\usepackage{amsmath}
\usepackage{bbm}
 \usepackage[foot]{amsaddr}

\usepackage{graphicx,xcolor}
\definecolor{green}{RGB}{89,169,58}
\definecolor{red}{RGB}{224,61,42}
\definecolor{blue}{RGB}{63,78,181}

\graphicspath{{fig/}}
\usepackage{subfigure,tikz}
\usetikzlibrary{arrows}

\newcommand{\N}{\ensuremath{\mathbb{N}}}
\newcommand{\R}{\ensuremath{\mathbb{R}}}
\newcommand{\Z}{\ensuremath{\mathbb{Z}}}
\newcommand{\E}{\ensuremath{\mathbb{E}}}
\renewcommand{\P}{\ensuremath{\mathbb{P}}}

\newcommand{\fref}[1]{Figure~\ref{#1}}

\newcommand{\ind}[1]{\ensuremath{\mathbbm{1}_{\left\{#1\right\}}}}
\newcommand{\diff}{\mathop{}\mathopen{}\mathrm{d}}
\newcommand{\cal}[1]{\ensuremath{\mathcal{#1}}}
\newcommand\croc[1]{\left\langle #1\right\rangle}
\newcommand\steq[1]{\stackrel{\text{\rm #1.}}{=}}

\def\eps{\varepsilon} 
\def\cadlag{c\`adl\`ag }

\newtheorem{proposition}{Proposition}
\newtheorem{definition}[proposition]{Definition}
\newtheorem{lemma}[proposition]{Lemma}
\newtheorem{remark}[proposition]{Remark}
\newtheorem{theorem}[proposition]{Theorem}

\numberwithin{proposition}{section}

\title[Noise induced stabilization in a CRN]{ Noise-induced stabilization in a chemical reaction network without boundary effects}

\author{Andrea Agazzi}
\email{andrea.agazzi@unibe.ch}
\author{Lucie Laurence}
\address{Institute of Mathematical Statistics and Actuarial Science, Department of Mathematics and Statistics, University of Bern, Alpeneggstrasse 22, 3012 Bern (CH)}
\email{lucie.laurence@unibe.ch}

\date{\today}

\begin{document}

\begin{abstract}
We present a chemical reaction network that is unstable under deterministic mass action kinetics, exhibiting finite-time blow-up of trajectories in the interior of the state space, but whose stochastic counterpart is positive recurrent. This provides an example of noise-induced stabilization of the model's dynamics arising due to noise perturbing transversally the divergent trajectories of the system that is independently of boundary effects.  The proof is based on a careful decomposition of the state space and the construction of suitable Lyapunov functions in each region.
\end{abstract}

\maketitle

\section{Introduction}\label{SecIntro}

Since the pioneering works of Feinberg \cite{Feinberg} and Horn and Jackson \cite{HornJackson}, a growing body of literature has focused on the relationship between the structure of a network of chemical reactions and its dynamics when modeled under the laws of mass action. These works can be naturally split in two groups based on the nature of the model used to describe such dynamics: deterministic or stochastic. In the former case, the dynamics is represented \emph{macroscopically} as the solution to a system of polynomial Ordinary Differential Equations (ODEs) modeling the time evolution of the concentration of the various species present in the system \cite{massactionbook}.
In the latter case, the \emph{microscopic} state of the network is represented by the number of molecules of each species at any given time, and the evolution of this state is modeled by a Continuous Time Markov Chain (CTMC), whose transition rates are determined by the structure of the network. These two lines of work were rigorously connected by the work of Kurtz \cite{Kurtz70,Kurtz71}, where the deterministic dynamics were identified as the Law of Large Number (LLN) limit of the stochastic ones (followed by Central Limit Theorem (CLT) \cite{Kurtz71,Kurtz78} and Large Deviation Principle (LDP) \cite{Agazzi1,Agazzi2,Agazzi3}) under an appropriate rescaling of the rates in a scaling parameter identified as the volume of the reactor.  Still, many questions about the relationship between these two models, further discussed in \cite{Anderson8}, remain open.

One of the most prominent application domains of the results outlined above is in systems biology, where one aims to model and understand the dynamics of living cells, interpreted as complex systems of chemical reactions, based on the structure of the associated reaction network. In this setting, however, a central role is played by the fact that reactions act on different timescales, allowing these systems to develop complex and hierarchical behavior that lies at the very core of their operational effectiveness. While the rigorous study of the multiscale behavior of such reaction networks is still in its infancy (see, e.g., \cite{LR23,LR24,LR24-2,Popovic06,Popovic14,Kumar17}), it is arguable that stochastic systems offer a critical advantage in modeling such multiscale behavior with respect to their deterministic counterparts. Indeed, stochastic models display a richer set of dynamical averaging properties in the long time limit, e.g., transitions between different dynamical attractors. Hence, a key step in the rigorous study of such multiscale systems is to guarantee sufficient stability, e.g., positive recurrence, of their dynamics on fast timescales, and the consequent existence of an invariant measure to which the process is expected to converge, on fast timescales.

The problem of the existence of an invariant measure has been extensively explored in the Stochastic Chemical Reaction Network Theory literature (see, e.g., \cite{Anderson1,Kang2,Cappelletti1} and references therein).
In these works, the main challenge results from the study of the dynamics near the boundary of state space, where it is known that the system behaves in a qualitatively different way depending on whether it is modeled deterministically or stochastically \cite{Anderson8,Cappelletti4,Agazzi4,Anderson9}.  The study of the asymptotic bulk dynamics, on the other hand, is typically carried out by treating the stochastic dynamics as a perturbation of the deterministic one and inferring stability of the stochastic model in the bulk from the stability of the deterministic limiting system in the spirit of the ODE method \cite{Borkar00}. 

\subsection*{Contributions} In this work, we provide a counterexample to this approach by presenting a network of chemical reactions that is \emph{unstable} (i.e., it has trajectories that diverge in finite time) in the bulk when modeled deterministically, but that becomes \emph{positive recurrent} when modeled stochastically. 
While a similar phenomenon was previously identified in \cite[Example 4]{Anderson9}, that result crucially depends on a carefully chosen set of reaction rate constants. The underlying stabilization mechanism in our construction, by contrast, is due to noise perturbing transversally the divergent trajectories of the system and as such does not hinge on parameter fine-tuning. Our proof
is based on a careful decomposition of the state space and the construction of
suitable Lyapunov functions in each region, adapting the approach developed in \cite{Mattingly2015} to our CTMC setting by combining it with weaker (non-pointwise) dissipation estimates in the spirit of \cite{Filonov1989}.
We further note that the dynamics of this reaction network appears to display intermittency (cf. Remark~\ref{r:excursions} below), a form of noise-amplifying chaotic behavior that is believed to be at the core of some important natural processes \cite{Alon,Bressloff}. 

 \subsection*{Setup and main result}

We consider a stochastic dynamical system modeling the evolution of a set of molecules of two species $S_1$ and $S_2$, that interact as prescribed by the following Chemical Reaction Network (CRN):
\begin{equation}\label{CRN}
	\begin{cases}
		3S_1+ 2S_2\rightharpoonup 5S_1+S_2, \quad  2S_1+ 3S_2\rightharpoonup S_1+5S_2, \\
		4S_1\rightharpoonup \emptyset, \quad 4S_2\rightharpoonup \emptyset\rightharpoonup S_1+S_2. 
	\end{cases}
\end{equation}
Concretely, denoting throughout by $\N$ the set of non negative integers, the state of this system at any time $t>0$ is captured by the vector $X(t) = (X_1(t), X_2(t)) \in \N^2$, counting the number of molecules of $S_1, S_2$ in the corresponding components. We model the evolution of $(X(t))$ as a continuous time Markov chain on $\N^2$. Under the laws of mass action \cite{massactionbook}, the $Q$-matrix $Q=(q(x,y), \ x,y\in \N)$ of $(X(t))$ is inherited from the network (\ref{CRN}), and its nonzero entries are given by
\begin{equation}\label{QMatrix}
	\begin{cases}
		q(x, x +2e_1-e_2)=x_1^{(3)}x_2^{(2)}\\
		q(x, x -e_1+2e_2)=x_1^{(2)}x_2^{(3)}\\
		q(x, x-4e_1)= x_1^{(4)}\\
		q(x, x-4e_2)= x_2^{(4)} \\
		q(x, x+e_1+e_2)=1
	\end{cases}
\end{equation}
Here, $(e_1, e_2)$ denote the orthonormal basis of $\R^2$, and for $y, p\in \N$, we define, throughout,
\begin{equation}
	y^{(p)}\steq{def}
	\begin{cases}
		\frac{y!}{(y-p)!} \quad \text{if}\quad y\geq p\\
		0 \quad \text{if} \quad y<p
	\end{cases}\,.
\end{equation}
While our results are independent on the choice of reaction rate constants, we set those constants to $1$ for all reactions for clarity of exposition.

\begin{figure}
\includegraphics[width=0.49\linewidth]{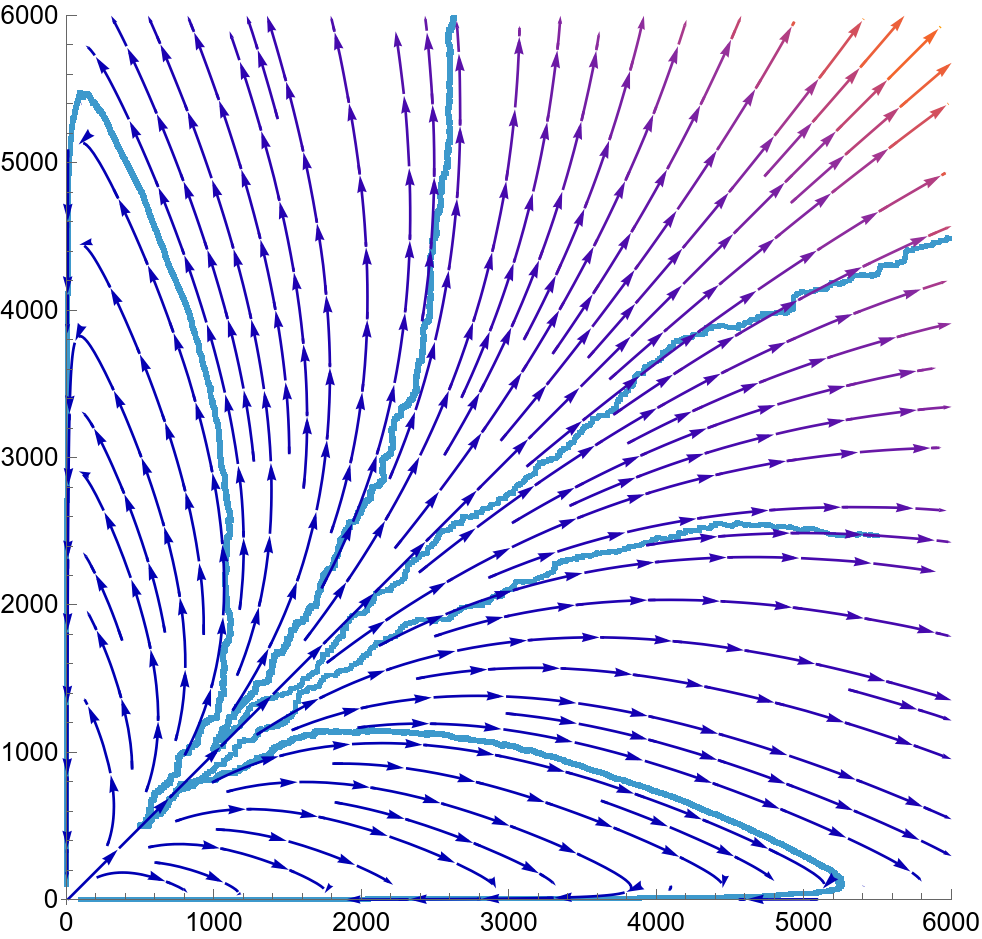}
\includegraphics[width=0.49\linewidth]{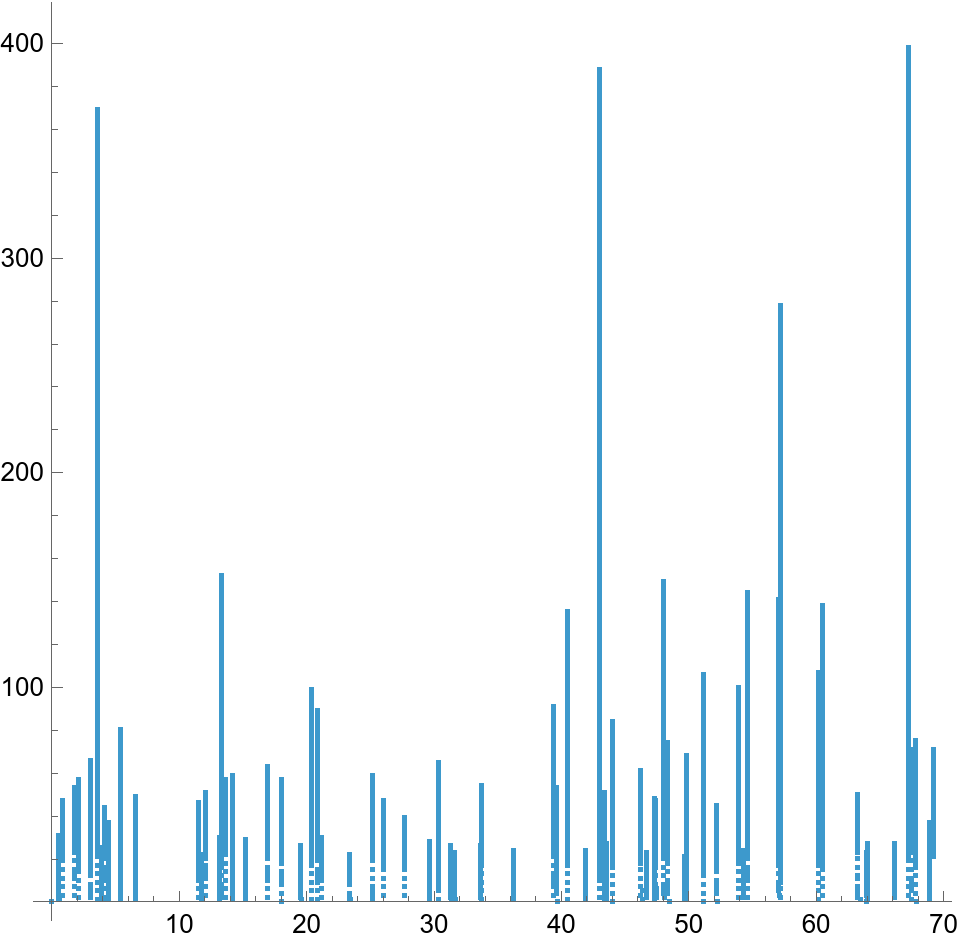}
\caption{Left: Flow lines of the vector field \ref{ODEs} and some sample paths of the stochastic process $X(t)$. Right: Sample path of $\|X(t)\|_1$ as a function of $t$, showing brusting intermittent behavior. All numerical experiments were run with Python.}
\label{f:vf}
\end{figure}

\begin{remark}\label{r:det}
Note that, when modeled deterministically, the system~\eqref{CRN} has trajectories that diverge in finite time while remaining arbitrarily far from the ``boundary'' of the state space $\{x_1 = 0\} \cup \{x_2 = 0\}$. Indeed, denoting throughout by $\mathbb R_{+}$ the set of nonnegative real numbers, under deterministic mass action kinetics the state $x(t) = (x_1(t), x_2(t)) \in \mathbb R_{+}^2$ of the network \eqref{CRN} evolves according to the system of ODEs
\begin{equation}\label{ODEs}
\begin{cases}
\frac d {dt} x_1 = 2 x_1^3 x_2^2 -  x_1^2 x_2^3 - 4 x_1^4+1\\ 
\frac d {dt} x_2 = - x_1^3 x_2^2 + 2 x_1^2 x_2^3 - 4 x_2^4+1
\end{cases}\,.
\end{equation}
The initial condition $x(0) = (x,x)$ yields a solution $x(t) = (f(x,t),f(x,t))$ where $f~:~\mathbb R \times \mathbb R_+ \to \mathbb R$ satisfies
$$\frac d {dt} f(x,t) = f(x,t)^5 - 4 f(x,t)^4+1\qquad f(x,0) = x$$
Hence, for $x > 4$, the trajectory $x(t)$ diverges in finite time. The flow induced by \eqref{ODEs} is depicted in \fref{f:vf}. Choosing different reaction rate constants preserves this divergent behavior, although on a different curve asymptotically satisfying $x_2 = \beta x_1$ for a $\beta > 0$.
\end{remark}

In contrast to the above remark and the fact that, whenever $X_1, X_2 \gg 1$, the dynamics of $X(t)$ is \emph{locally} well approximated by the trajectories of \eqref{ODEs} with the same initial condition, our main result states that the long-time stability properties of the process $(X(t))$ are different from the ones of the deterministic model:
\begin{theorem}\label{StabThm}
	The process $(X(t))$ associated to  the CRN \eqref{CRN} is positive recurrent. 
\end{theorem}
In particular, since the process $(X(t))$ is manifestly irreducible, it possesses a unique invariant measure. 

\begin{remark}\label{r:excursions} As can be observed from \fref{f:vf}, the stability of this system is a consequence of the perturbative instability of the divergent trajectories of the deterministic system identified in Remark~\ref{r:det}. In particular, when initialized sufficiently close to the unstablle manifold $\{(x_1,x_2) \in \N^2~:~x_1 = x_2\}$, the process $(X(t))$ remains locally close to the the solution of \eqref{ODEs}, but is then pushed by the random fluctuations onto a stabilizing trajectory, which it approximately follows until it returns to a compact.
As a consequence of this effect, this system displays intermittent bursting behavior: macroscopic excursions appear to occur as a Poisson process, and the tail distribution of their excursion length can be predicted through an estimate similar to \cite{Mattingly2015,Herzog15}.
\end{remark}

The remainder of the paper is structured as follows. In Section~\ref{SecModel} we outline the strategy of the proof of Theorem~\ref{StabThm} defining a partition of the state space into regions where the behavior of the process is qualitatively different. Then, in Section~\ref{SecCone} we obtain the desired estimates in the neighborhood of the unstable trajectory identified in Remark~\ref{r:det} and in Section~\ref{SecSides} we do the same on the complement of this neighbborhood in $\N^2$. 	
\section{Structure of the proof}\label{SecModel}

We prove Theorem~\ref{StabThm} by identifying a function that is dissipated, on average, by the trajectories of our system in a sufficiently weak sense.  In Section~\ref{s:filonov} we define the dissipation condition that we aim to establish, introduced by Filonov \cite{Filonov1989}. In Sections~\ref{s:EnergyFunction} and~\ref{s:ProofThSketch} we partition the state space asymptotically into regions where the dynamics is dominated by different leading order terms, and define a candidate energy function on each element of such partition. Finally, we state the partial dissipation results on each such region (proven in the following sections) which, combined, prove our main Theorem~\ref{StabThm}. Our construction is inspired by the one presented in \cite{Mattingly2015}, combining it with non-pointwise dissipation estimates allowing to significantly streamline the proof in the present setting.

\subsection{Filonov's criterion}\label{s:filonov}

It is easy to check that the process is irreducible on $\N^2$. 
The proof of Theorem~\ref{StabThm} is carried out using Filonov's criterion, see~\cite{Filonov1989} or~\cite{LR23}: 
\begin{proposition}[Filonov's criterion]{}\label{p:filonov}
	If there exists a function $V: \N^2\to \R_+$ with 
	\begin{equation}\label{CondV}
		\lim_{\|x\|_1\to +\infty} V(x)=+\infty, 
	\end{equation}
	an integrable stopping time $\tau$, and some constants $\gamma, K>0$ such that 
	\begin{equation}\label{InFilo}
		\E_x\left(V(X(\tau))\right)-V(x)\leq -\gamma\E_x(\tau) \qquad \text{for all } \|x\|_1\geq K\,,
	\end{equation}
	then the process $(X(t))$ is positive recurrent. 
\end{proposition}
A function $V$ that verifies Condition~\eqref{CondV} is called an \emph{energy function}. 
Filonov's criterion is an extension of the Lyapunov theorem. 
 Indeed, consider a real-valued function $(x(t))$ on $\R^2$ that is \cadlag (right continuous and it has left-limits) everywhere on $\R_+$, for $t>0$, ans let $x(t{-})$ denote the left limit at time $t$. Then, if one chooses $\tau$ as the first jumping time $t_1$ of the process, i.e., 
\begin{equation}\label{eq:t1}
	t_1\steq{def} \inf\{t\geq 0: \quad X(t)\neq X(t{-})\}, 
\end{equation}
one gets for $x\in \N^2$
\begin{align*}
	\E_x\left(V(X(t_1))-V(x)\right)=\E_x\left(\int_0^{t_1}\cal{L}(V)(X(s{-}))\diff s\right)=\cal{L}(V)(x)\E_x\left(t_1\right), 
\end{align*}
where $\cal{L}$ is the infinitesimal generator of the process $(X(t))$. Therefore, when $\tau=t_1$, Relation~\eqref{InFilo} is equivalent to 
\begin{equation}\label{InLyap}
	\cal{L}(V)(x)\leq -\gamma. 
\end{equation}

To show that the process $(X(t))$ verifies Filonov's criterion, one strategy is to set $\tau=t_1$ at each state for which Relation~\eqref{InLyap} is verified and choose a more complicated stopping time only for the initial states where it is necessary. 
Note that the process is symmetric in both coordinates, so it is sufficient to show the criterion for states $x\in \N^2$ such that $x_1\geq x_2$.

\subsection{Definition of the energy function}\label{s:EnergyFunction}

Depending on the order of magnitude of the number of molecules of each species, different reactions can dominate the evolution of the system. To take this feature into account, we partition the state space into three subspaces: a cone $\cal{C}$ centered on the diagonal, and the left and right complements of this set in $\N^2$, respectively denoted by $\cal{S}_L$ and  $\cal{S}_R$.

\begin{definition}
	For some constant $p\in(0,1)$ specified later,  we set 
	\begin{equation}\label{DefSets}
		\begin{cases}
			\cal{C}&=\left\{(x_1, x_2)\in \N^2: px_2 < x_1<\frac{x_2}{p} \right\}, \\
			\cal{S}_R&=\left\{(x_1, x_2)\in \N^2: x_2\leq px_1\right\}, \\
			\cal{S}_L&=\left\{(x_1, x_2)\in \N^2:  x_1\leq p x_2 \right\}.  
		\end{cases}
	\end{equation}
\end{definition}

The sets are chosen symmetric with respect to the line $x_1=x_2$ consistently with the symmetry of the system. For each of these subsets, we can identify the dominating reactions. In the cone, reactions $2S_1+3S_2\rightharpoonup S_1+5S_2$ and $3S_1+2S_2\rightharpoonup 5S_1+S_2$ have the larger rates, whereas in $\cal{S}_L$, it is either reaction $2S_1+3S_2\rightharpoonup S_1+5S_2$ or reaction $4S_2\rightharpoonup \emptyset$ that dominates, respectively in $\cal{S}_R$, either reaction $3S_1+2S_2\rightharpoonup 5S_1+S_2$ or $4S_1\rightharpoonup \emptyset$ dominates. 

We now define a simple energy function, for which we aim to prove that Relation~\eqref{InFilo} holds. We define this function piecewise on the partition given above, and such that in each subset, the dominating reactions are dissipative. 
%
%
%
%
 Denoting throughout for $x,y\in \R$, $x\wedge y= \min\{x,y\}$ and $x\vee y= \max\{x, y\}$, we set
\begin{equation}\label{DefV}
	V:(x_1, x_2)\mapsto 
	\begin{cases}
		V_C(x) \quad \text{when} \quad x\in\cal{C}, \\
		V_R(x) \quad \text{when} \quad x\in \cal{S}_R,\\
		V_L(x) \quad \text{when} \quad x\in \cal{S}_L, 
	\end{cases}
\end{equation}
where for $x\in \R_+^2$, 
\begin{align*}
	V_R(x)\steq{def}c_S(x_1+6x_2), \quad V_L(x)\steq{def}c_S(6x_1+x_2),
\end{align*}
and
\begin{align*}
 V_C(x)\steq{def}\frac{(x_1+x_2)^3}{(|x_1-x_2|\vee \beta \sqrt{x_1+x_2})^2}\,,
\end{align*}
for a  $\beta>0$ that will be fixed later and $c_S$ such that $V_L(x)=V_C(x)$, respectively $V_R(x)=V_C(x)$ at the interface $\{x_2 = x_1/p\}$, respectively $\{x_2 = p x_1\}$: 
\[
	c_S\steq{def} \frac{(1+p)^3}{(1-p)^2(1+6p)}\,. 
\]

The rest of the article is dedicated to the proof of Relation~\eqref{InFilo} with the energy function defined in Relation~\eqref{DefV}.

\subsection{Proof of Theorem~\ref{StabThm}}\label{s:ProofThSketch}

To prove our main result we further partition the cone $\mathcal C$ into sets where the dynamics are, asymptotically, in a diffusive or in a transport regime, respectively due the leading order reactions in this set having comparable intensities or one being stronger, in an appropriate sense, than the other. 
\subsubsection*{Partitioning the cone $\mathcal C$}

To study the process in the cone $\cal{C}$, we consider the following convenient change of coordinates.
\begin{definition}
	If $x=(x_1, x_2)\in \N^2$, we introduce 
	\[
		\Phi(x_1, x_2)=(r(x_1, x_2), d(x_1, x_2))=(x_1+x_2, x_1-x_2). 
	\]
	We call these coordinates the \emph{axial coordinates}.
	The inverse function $\Phi^{-1}$ is defined on 
	\[
		\Phi(\N^2)=\{(r,d)\in \N\times \Z: \quad |d|\leq r\}. 
	\]
	For a set $A\subset \N^2$, we set $\tilde{A}=\Phi(A)$. 
\end{definition}

We introduce the process $(R(t), D(t))$ on $\N\times \Z$ as	
\begin{equation}\label{DefRD}
	(R(t), D(t))\steq{def}(\Phi(X_1(t), X_2(t))).
\end{equation}
The process $(R(t), D(t))$ is a CTMC, whose transitions can be deduced from the $Q$-matrix defined in Relation~\eqref{QMatrix}. Its infinitesimal generator is given by 
\begin{equation}\label{InfGenAxial}
	\cal{A}(f)(r, d)=\cal{L}(f\circ\Phi)(\Phi^{-1}(r,d)). 
\end{equation}
for a real function $f$ defined on $\N\times \Z$, and $(r,d)\in \Phi(\N^2)$. 
The process leaves the cone $\cal{C}$ when $|D(t)|\geq  q R(t)$ for $q=(1-p)/(1+p)>0$, where $p$ is given in Relation~\eqref{DefSets}.

The inequalities~\eqref{InFilo} and~\eqref{InLyap} are easily translated into the axial coordinates, using the function $V\circ\Phi^{-1}$ as a Lyapunov function, i.e., 
\[
	\tilde{V}_C(r,d)\steq{def} V_C\circ \Phi^{-1}(r,d)=\frac{r^3}{(d\vee\beta \sqrt{r})^2}. 
\]

In Section~\ref{SecSides}, we will see that Relation~\eqref{InLyap} holds for states on the side of the cone. More precisely, the  states $(r,d)\in \tilde{\cal{C}}$ such that $|d|\geq \eta_0\sqrt{r}$ for some $\eta_0$ large enough, specified later. 

We therefore define a partition of $\tilde{\cal{C}}$ as follow: 
\begin{equation}\label{DefIntCone}
	\begin{cases}
		\tilde{\cal{C}}_{c}\steq{def}\{(r, d)\in \tilde{\cal{C}}: |d|\leq \eta_0 \sqrt{r}\}=\tilde{\cal{C}}_c(\eta_0),\\
		\tilde{\cal{C}}_R\steq{def}\{(r, d)\in \tilde{\cal{C}}: d>\eta_0 \sqrt{r}\},\\
		\tilde{\cal{C}}_L\steq{def}\{(r, d)\in \tilde{\cal{C}}: d<-\eta_0 \sqrt{r}\}, 
	\end{cases}
\end{equation}
and their equivalent in the Cartesian coordinates: $\cal{C}_c=\Phi^{-1}(\tilde{\cal{C}}_c)$, $\cal{C}_L=\Phi^{-1}(\tilde{\cal{C}}_L)$ and  $\cal{C}_R=\Phi^{-1}(\tilde{\cal{C}}_R)$. For $\eta>0$, we also set $\tilde{\cal{C}}_c(\eta)$ replacing $\eta_0$ by $\eta$ in the definition.

\subsubsection*{The transport regime}
In the complement of the set $\cal{C}_c$ it is sufficient to choose $\tau = t_1$ from Relation~\eqref{eq:t1}, as discussed after the statement of the main theorem. Here and throughout, for $y\in \N^2$ we define $\|y\|_1=y_1+y_2$.

\begin{proposition}\label{p:filonov2}
	The exist some $\gamma>0$, $N_0\geq 1$ such that for all $x_N\in \mathbb N^2\setminus \cal{C}_c$ with $\|x_N\|_1=N\geq N_0$, the following inequality holds :
	\begin{equation}
		\E_{x_N}\left(V(X(t_1))\right)-V(x)\leq -\gamma\E_{x_N}(t_1), 
	\end{equation}
	where $V$ is defined in Relation~\eqref{DefV} and $t_1$ in Relation~\eqref{eq:t1}.  
\end{proposition}

\begin{proof}[Proof of Proposition~\ref{p:filonov2}] This follows from combining Lemma~\ref{PropSide1} and Lemma~\ref{PropSide2}, proven in Section~\ref{SecSides}.
\end{proof}

\subsubsection*{The diffusive regime} It remains to define the stopping time $\tau$ in Relation~\eqref{InLyap} for the states in $\tilde{\cal{C}}_c$. This is the hardest part of the proof. In Section~\ref{SecCone}, we show that the energy of the system decreases if the process deviates from the axis $d=0$ (or $x_1=x_2$ in the Cartesian coordinates) and that here, it is the randomness of the process that helps it get away from the axis. This is summarized in the following proposition.

\begin{proposition}\label{p:filonov1}
	There exist some $\gamma>0$, $N_0\geq 1$ and an integrable stopping time $\tau$ such that for all $x_N\in \cal{C}_c$ with $\|x_N\|_1=N\geq N_0$, the following inequality holds :
	\begin{equation}\label{EqFiloGronwall}
		\E_{x_N}\left(V(X(\tau))\right)-V(x)\leq -\gamma\E_{x_N}(\tau), 
	\end{equation}
	where $V$ is defined in Relation~\eqref{DefV}.  
\end{proposition}

This result is proven is Section~\ref{SecCone}. These two propositions allow to prove our main result:

\begin{proof}[Proof of Theorem~\ref{StabThm}] Proposition~\ref{p:filonov1} and Proposition~\ref{p:filonov2} jointly show that the conditions of Proposition~\ref{p:filonov} hold and therefore prove our main result Theorem~\ref{StabThm}.\end{proof}
{}

\section{Interior of the cone}\label{SecCone}

In this section, we prove Proposition~\ref{p:filonov1}, i.e., that Relation~\eqref{InFilo} holds for the states in $\cal{C}_c$. Axial coordinates will be used throughout. 
A convenient stopping time to obtain the desired inequality can be defined (up to technical adjustments such as further stopping and process approximation) as 
\begin{equation}
	\inf\{t\geq 0: |D(t)|\geq \eta_1 \sqrt{R(t)}\}, 
\end{equation}
for some $\eta_1>\eta_0$, i.e., as the exit time from a larger interior cone $\cal{C}_c(\eta_1)$.

The proof is carried out in several steps: 
\begin{itemize}
	\item First, we apply a time change to $(X(t))$, to capture the timescale of the fastest reactions, obtaining a slowed down version $(Z(t))$ that evolves on a timescale of order 1. 
	\item Then we introduce an auxiliary process $(Y(t))$, approximating the dynamics of $(Z(t))$ to leading order, by only taking into account the dominant transitions in the interior cone $\cal{C}_c$ (i.e., the ones with dominant rates). We show some limit theorems on the sample paths of $(Y(t))$. 
	\item These limit theorems allows us to estimate the value of the energy of $(Y(t))$ at a stopping time $\tau_Y$ carefully chosen, and conclude to Relation~\eqref{InFilo} for $(Y(t))$. 
	\item We show that the processes $(Y(t))$ and $(Z(t))$ have a similar behavior, by showing that they stay sufficiently close on the timescale studied. 
	\item To conclude, we only check that the time change does affect the relevant estimate Relation~\eqref{InFilo}. 
\end{itemize}

Before proceeding with the proof, we introduce some necessary notation.  We set $\croc{-,-}$ the Euclidean product in $\R^2$. For a Poisson process $\cal{P}$ on $\R_+^2$ with intensity measure the Lebesgues measure on $\R_+^2$, we will use the notation 
\[
	\cal{P}(A, \diff t)=\int_{x\in A} \cal{P}(\diff x, \diff t),  
\]
for $A\in \cal{B}(\R_+)$. 

\subsection{Definition of the auxiliary process}

Let $N\geq 1$, and consider the process $(X_N(t))$ starting from an initial state $x_N\in \cal{C}_c(\eta_0)$ with $\|x_N\|_1=N$, and whose $Q$-matrix is given by Relation~\eqref{QMatrix}. 

We introduce the process $(Z^N(t))$ starting at $x_N$, whose $Q$-matrix is given by 
\begin{equation}\label{QMatrixZ}
	\begin{cases}
		q_{Z}(z, z +2e_1-e_2)=z_1-2, \\
		q_{Z}(z, z -e_1+2e_2)=z_2-2, \\
		q_{Z}(z, z -4e_1)=z_1^{(4)}/z_1^{(2)}z_2^{(2)},\\
		q_{Z}(z, z -4e_2)=z_2^{(4)}/z_1^{(2)}z_2^{(2)},\\
		q_{Z}(z, z +e_1+e_2)=1/z_1^{(2)}z_2^{(2)}. 
	\end{cases}
\end{equation} 
The process $(Z(t))$ is not well defined when it reaches ``small states'', i.e., when $z_1\leq 2$ or $z_2\leq 2$. However, we will only study it in $\cal{C}$, where both coordinates are large.

Define for $t\geq 0$, the stopping time 
\begin{equation}\label{EqLambda}
	\lambda^N_t=\inf\left\{s>0: \int_0^s \frac{1}{Z^N_1(u)^{(2)}Z^N_2(u)^{(2)}}\diff u\geq t\right\}, 
\end{equation}
and let $(\tilde{X}^N(t))=(Z^N(\lambda_t^N))$, then one can check that $(\tilde{X}^N(t))$ and $(X^N(t))$ have the same distribution (see \cite{Rogers1}, Section III.21). 

In the cone $\cal{C}_c(\eta_1)$, two transitions, associated to the reactions $3A+2B\rightharpoonup 5A+B$ and $2A+3B\rightharpoonup A+5B$, dominate the evolution of the process $(Z^N(t))$. 
We define the continuous time Markov chain $(Y^N(t))$, simplifying the process $(Z^N(t))$ by only taking into account these two leading order transitions. 
More specifically, we set $Y^N(0)=x_N$ and define $(Y^N(t))$ as the solution of the stochastic differential equation (SDE)
\begin{equation}\label{SDEY}
	\begin{cases}
		\diff Y_1(t)= 2\cal{P}_1((0, Y_1(t{-})-2), \diff t)-\cal{P}_2((0, Y_2(t{-})-2), \diff t),\\
		\diff Y_2(t)= -\cal{P}_1((0, Y_1(t{-})-2), \diff t)+2\cal{P}_2((0, Y_2(t{-})-2), \diff t).
	\end{cases}
\end{equation}
where $\cal{P}_1$ and $\cal{P}_2$ are two independent Poisson processes on $\R_+^2$ with intensity measure the Lebesgue measure on $\R_+^2$. We will show in Section~\ref{SecBack} that when $N$ gets large, $(Y^N(t))$ approximates well the evolution of the process $(Z^N(t))$. 
The two next sections focus on the process $(Y^N(t))$. 

\subsection{Scaling results}

We introduce throughout the processes in the axial coordinates as follow: 
\begin{equation}\label{DefRDY}
	(R_Y(t), D_Y(t))\steq{def} \Phi(Y(t)). 
\end{equation}

\begin{definition}
	We introduce $(L(t))$ the divergent Ornstein-Uhlenbeck process, starting at $d_0\geq 0$ and solving the SDE 
	\begin{equation}\label{EqOU}
		\diff L(t)= \frac{5}{2} L(t)\diff t + 3 \diff B(t)
	\end{equation}
	with $(B(t))$ a Brownian process.  	
\end{definition}

The following result controls the trajectory of $(R_Y(t), D_Y(t))$ in the large $N$ regime.
\begin{proposition}\label{PropScaling}
	Let $N\geq 1$, $d_N\in \Z$ such that $(N, d_N)\in \tilde{\cal{C}}_c$, $(R^N_Y(t), D_Y^N(t))$ the process defined in Relation~\eqref{DefRDY}, with initial condition $(R^N_Y(0), D_Y^N(0))=(N, d_N)$. 
	
	For $T>0$ and $k\geq 1$, the convergence in distribution
	\[
		\lim_{N\to +\infty} \left( \left(\frac{R_Y^N(t)}{N}\right)^k, t\in [0,T]\right) =(e^{kt}, \ t\in [0, T])
	\]
	holds, and if $\lim_{N\to +\infty} (d_N/\sqrt{N})= d_0$,  the convergence in distribution
	\[
		\lim_{N\to +\infty} \left( \frac{D_Y^N(t)}{\sqrt{R_Y^N(t)}}, t\in[0, T]\right) =(L(t), \ t\in [0,T])
	\]
	holds.
\end{proposition}

\begin{proof}
	\emph{Convergence of $(R_Y^N(t)/N)$:}\\
	Let $T>0$. Integrating Relation~\eqref{SDEY}, we have 
	\begin{equation}\label{Rel1}
		R_Y^N(t)= N+ \int_0^t (Y_1^N(s)+Y_2^N(s)-4)\diff s + M_R^N(t), 
	\end{equation}
	where $(M_R^N(t))$ is a martingale whose predictable quadratic variation is given by 
	\[
		\croc{M_R^N}(t)=\int_0^t(R_Y^N(s)-4)\diff s. 
	\]
	Therefore, introducing for $N\geq 1$ the scaled process $(\overline{R}_Y^N(t))=(R_Y^N(t)/N)$, we have 
	\begin{equation}
		\E_{x_N}\left(\overline{R}_Y^N(t)\right)\leq 1+\int_0^t \E_{x_N}\left(\overline{R}_Y^N(s)\right). 
	\end{equation}
	Using Gronwall's Lemma, this leads to 
	\[
		\sup_{N\geq 1}\sup_{t\leq T}\E\left(\overline{R}_Y^N(t)\right)\leq e^T. 
	\]
	Therefore, using Doob's inequality, the relation
	\[
		\lim_{N\to+\infty}\E\left(\sup_{t\leq T}\frac{|M_R^N(t)|}{N}\right)=0
	\]
	holds. 
	Using Relation~\eqref{Rel1} and Gronwall's Lemma, we get 
	\[
		\sup_{N\geq 1}\E_{(N, d_N)}\left(\sup_{t\leq T} \overline{R}_Y^N(t)\right)<2e^T. 
	\]
	We now define the modulus of continuity of $\overline{R}_Y^N(t)$, see~\cite{Billingsley}: for $\delta>0$, 
	\[
		\omega^N_R(\delta)\steq{def}\sup\left\{| \overline{R}_Y^N(t)-\overline{R}_Y^N(s)|, \ s,t\in [0,T], \ |s-t|\leq \delta\right\}. 
	\]
	Relation~\eqref{Rel1} gives
	\[
		\omega^N_R(\delta)\leq \frac{4\delta}{N}+ \delta\sup_{t\leq T}\overline{R}_Y^N(t) + 2\sup_{t\leq T}|M_R^N(t)|, 
	\]
	and therefore for $\eps>0$, and $\eta>0$, one can find $N_0\geq 1$ and $\delta_0>0$ such that for $\delta<\delta_0$ and $N\geq N_0$, 
	\[
		\P\left(\omega^N_R(\delta)>\eta\right)\leq \eps, 
	\]
	and $(\overline{R}_Y^N(t))$ is tight. Relation~\eqref{Rel1} shows that the limit $(r(t))$ verifies $\dot{r}(t)=r(t)$, which concludes the proof of the first convergence for $k=1$. Extension to $k\geq 2$ is straightforward.  
	
	\vspace{0.5cm}
	
	\emph{Convergence of $\left(D_Y^N(t)/\sqrt{R_Y^N(t)}\right)$:}\\
	We set $g: (r, d)\mapsto d/\sqrt{r}$. 
	The process 
	\[
		(L_N(t))\steq{def}\left(\frac{D_Y^N(t)}{\sqrt{R_Y^N(t)}}\right)
	\]
	verifies the Relation 
	\begin{multline}\label{EqLN}
		L_N(t)=\frac{d_N}{\sqrt{N}} +\int_0^t\frac{R_Y^N(s){+}D_Y^N(s)-4}{2}\left(g(R_Y^N(s)+1, D_Y^N(s)+3)-g(R_Y^N(s), D_Y^N(s)) \right)\diff s\\
		\hfill +\int_0^t \frac{R_Y^N(s){-}D_Y^N(s){-}4}{2} \left(g(R_Y^N(s)+1, D_Y^N(s)-3)-g(R_Y^N(s), D_Y^N(s))\right)\diff s +M_L^N(t), 
	\end{multline}
	where $(M_L^N(t))$ is a square integrable martingale whose increasing process is 
	\begin{multline}\label{CrocML}
		\croc{M_L^N}(t)=\int_0^t \frac{R_Y^N(s)+D_Y^N(s)-4}{2} \left(g(R_Y^N(s)+1, D_Y^N(s)+3)-g(R_Y^N(s), D_Y^N(s)) \right)^2\diff s\\
		+\int_0^t\frac{R_Y^N(s)-D_Y^N(s)-4}{2} \left(g(R_Y^N(s)+1, D_Y^N(s)-3)-g(R_Y^N(s), D_Y^N(s))\right)^2\diff s.  
	\end{multline}
	
	Using a Taylor expansion, one can show that for $a\in \{-3,3\}$, there exists some $r^*\geq 1$ such that,  if $r\geq r^*$ and $d\geq 1$, 
	\begin{equation}\label{ApproxLN}
		\left|g(r+1, d+a)-g(r,d)-\frac{a}{\sqrt{r}}+\frac{d+a}{2r^{3/2}}\right|\leq \frac{C_g d}{r^{5/2}}
	\end{equation}
	for some constant $C_g>0$. 
	
	Using this inequality in Relation~\eqref{CrocML}, one can show that
	\[
		\lim_{N\to +\infty} \left(\croc{M_L^N}(t), t\in[0,T]\right)=(9t, \ t\in [0, T]), 
	\]
	in distribution, and deduce, using~\cite{KurtzEthier}, Theorem 1.4, that $(M_L^N(t))$ converges in distribution to $(B(9t))$, where $(B(t))$ is the standard Brownian motion.

	Making the approximation 
	$$
		g(r+1, d+a)-g(r,d)\approx \frac{a}{\sqrt{r}}-\frac{d+a}{2r^{3/2}}
	$$
	in Relation~\eqref{EqLN}, one gets 
	\[
		L_N(t)=\frac{d_N}{\sqrt{N}} +\int_0^t\frac{5}{2}L_N(s)\diff s+\Delta_L^N(t) +M_L^N(t), 
	\]
	where $\Delta_L^N(t)$ contains the remaining terms after the approximation is made. 
	
	Using the inequality~\eqref{ApproxLN}, we get
	\begin{align*}
		|\Delta_L^N(t)|&\leq \int_0^t\frac{R_Y^N(s){+}D_Y^N(s)-4}{2}\frac{C_gD_Y^N(s)}{(R_Y^N(s))^{5/2}}\diff s\\
					&\quad +\int_0^t\frac{R_Y^N(s){-}D_Y^N(s)-4}{2}\frac{C_gD_Y^N(s)}{(R_Y^N(s))^{5/2}}\diff s\\
					&\leq\int_0^t\frac{C_gD_Y^N(s)}{(R_Y^N(s))^{3/2}}\diff s. 
	\end{align*}
	
		Using the same steps as the proof of the convergence of $(\overline{R}^N_Y(t))$, we can show that $(D_Y^N(t)/N^{3/2})$ converges in distribution to $0$. With the convergence of $(R_Y^N(t)/N)$, it is therefore straightforward that $(\Delta_L^N(t))$ converges in distribution to $(0)$. 
		
	The tightness of the process $(L_N(t))$ can now be proved using its modulus of continuity, and similar arguments as for the tightness of $(R_Y^N(t))$. If $(L(t))$ is a limiting point of this sequence, $(L(t))$ is continuous, see~\cite{Billingsley}, and is the unique solution of, for all $t\in [0, T]$ 
	\[
		\diff L(t)= \frac{5}{2}L(t) \diff t + 3\diff W(t), \quad L(0)=d_0. 
	\]
	where $(W(t))$ is a standard Brownian motion, see~\cite{Rogers2} for example. 
	This concludes the proof of the proposition. 
\end{proof}

\subsection{Filonov criterion for the process $(Y^N(t))$ in the cone}

\begin{definition}
	We introduce the two following stopping times: \\
	For some $\eta_1>0$ specified later, $N\geq 1$ and $(R_Y^N(0),D_Y^N(0))=(N, d_N)\in \tilde{\cal{C}_c}$, we set 
	\begin{equation}\label{DefTauY}
		\tau_Y^N\steq{def}\inf\{ t\geq 0:\ |D_Y^N(t)|\geq \eta_1\sqrt{R_Y^N(t)}\}=\inf\{t\geq 0: |L_N(t)|\geq \eta_1\}, 
	\end{equation}
	and 
	\[
		\tau_L\steq{def}\inf\{ t\geq 0:\  |L(t)|\geq \eta_1\}. 
	\]
\end{definition}

In the next proposition, adapted from \cite{Mattingly2015}, we bound the tails of the stopping time $\tau_L$ from above: 
\begin{proposition}\label{PropStopTime}
	There exists a constant $C_\tau$ independent of $\eta_1$ such that  
	\[
		\E\left(e^{2\tau_L}| L(0)=0\right)< 2 +C_\tau \eta_1. 
	\]	
\end{proposition}

\begin{proof}
	Note that if $(L(t))$ starts at $d_0\in [-\eta_1, \eta_1]$ and is solution of~\eqref{EqOU}, then for $t\geq 0$, $L(t)$ is normally distributed with mean $d_0e^{5t/2}$ and variance ${9(e^{5t}-1)/5}$. Therefore, for $s>0$, 
	\begin{equation}\label{MajQueueTau}
		\P_0\left(\tau_L>s\right)\leq \P_0\left(\sup\{|L(v)|: \ v\leq s\}<\eta_1\right)
		\leq \P_0(L(s)< \eta_1)\leq \frac{\eta_1K_0}{\sqrt{e^{5s}-1}},
	\end{equation}
	where $K_0\steq{def}2\sqrt{5}/3\sqrt{2\pi}$, and in the last inequality we used the explicit expression of the distribution of $L(t)$.
	This leads to
	\begin{align*}
		\E_{0}(e^{2\tau_L})&=\int_0^{+\infty} \P_{0}(e^{2\tau}>u)\diff u \\
			&\leq 2+ \int_2^{+\infty} \frac{\eta_1 K_0}{\sqrt{u^{5/2}-1}}\diff u\\
			&\leq 2+\eta_1 C_\tau, 
	\end{align*}
	 for some $C_\tau$ independent of $\eta_1$. 
\end{proof}

\begin{proposition}\label{PropStabY2}
	For some choice of $0<\eta_0<\eta_1$, there exists $T>0$, $\gamma>0$ and $N_0\geq 1$ such that for $N\geq N_0$, for all $d_N\in \Z$ such that $(N, d_N)\in \tilde{\cal{C}}_c$, 
	\begin{equation}\label{EqPropStabY}
		\E_{(N, d_N)}\left(\tilde{V}(R(\tau_Y^N\wedge T), D(\tau_Y^N\wedge T))-\tilde{V}(N, d_N)\right)\leq -\gamma. 
	\end{equation}
\end{proposition}

\begin{proof}
	Let $\eps_0>0$. Recall that $\tilde{V}(r,d)=r^3/(d\vee \beta\sqrt{r})^2$. 
	
	We have: 
	\begin{align*}
		N^{-2}\E_{(N, d_N)}&\left(\tilde{V}_C(R_Y^N(\tau_Y^N\wedge T),D_Y^N(\tau_Y^N\wedge T))\right)\\
		&=\E_{(N, d_N)}\left((\overline{R}_Y^N(\tau_Y^N\wedge T))^2(L_N(\tau_Y^N\wedge T)\vee \beta)^{-2}\right)\\
		&\leq \beta^{-2}\E_{(N, d_N)}\left(\sup_{t<T}\left|\overline{R}_Y^N(t)^2-e^{2t}\right|\right)
		+\E_{(N, d_N)}\left(e^{2(\tau_Y^N\wedge T)}(L_N(\tau_Y^N\wedge T)\vee \beta)^{-2}\right)\\
		&\leq \beta^{-2}\E_{(N, d_N)}\left(\sup_{t<T}\left|\overline{R}_Y^N(t)^2-e^{2t}\right|\right)
		+\frac{e^{2T}}{\beta^2}\P_{(N, d_N)}\left(\tau_Y^N>T\right)\\
		&\quad + \E_{(N, d_N)}\left(e^{2(\tau_Y^N\wedge T)}\eta_1^{-2}\right). 
	\end{align*} 
	We deal with each term separately, choosing the constants $\eps_0>0$ small enough, $T>0$ and $\eta_1>\eta_0$ large enough so that: 
	\begin{equation}\label{CondConst}
		\eps_0+\frac{2+C_\tau\eta_1}{\eta_1^2}<\frac{1}{\eta_0^2} \quad \text{and}\quad \frac{\eta_1e^{2T}K_0}{\beta^2\sqrt{e^{5T}-1}}<\eps_0/4, 
	\end{equation}
	where $C_\tau$ is defined in Proposition~\ref{PropStopTime} and $K_0$ in Relation~\eqref{MajQueueTau}.

	For the second term, we have:
	\begin{align*}
		\frac{e^{2T}}{\beta^2}\P_{(N, d_N)}&\left(\tau_Y^N>T\right)\leq \frac{e^{2T}}{\beta^2}\max_{0\leq k\leq 10}\P_{(N, k)}\left(\tau_Y^N>T\right)\\
		&\leq \frac{e^{2T}}{\beta^2}\left|\max_{0\leq k\leq 10}\P_{(N, k)}\left(\tau_Y^N>T\right)- \P_0\left(\tau_L>T\right)\right|+\frac{e^{2T}}{\beta^2}\P_0\left(\tau_L>T\right)\\
	&\leq \eps_0/4+\frac{\eta_1e^{2T}K_0}{\beta^2\sqrt{e^{5T}-1}}.
	\end{align*}
	The first inequality is the consequence of some monotonicity of the stopping time $\tau_Y^N$ according to the initial state proven in Lemma~\ref{LemmaOrderTau}, see below. The second inequality, valid for $N\geq N_{02}$ for some $N_{02}$ large enough, comes from combining Relation~\eqref{MajQueueTau} and the convergence of $\tau_Y^N$ to $\tau_L$, consequence of the convergence in distribution of $(L_N(t))$ to $(L(t))$. 
	
	For the first term, Proposition~\ref{PropScaling} allows us to find some $N_{01}\geq 1$ such that 
	\begin{align*}
		\sup_{d_N\in[-\eta_0\sqrt{N}, \eta_0\sqrt{N}]}\beta^{-2}\E_{(N, d_N)}\left(\sup_{t<T}\left|\overline{R}_Y^N(t)^2-e^{2t}\right|\right)<\eps_0/4. 
	\end{align*}
	
	For the third term finally, one can find $N_{03}\geq 1$ such that 
	\begin{align*}
		\frac{1}{\eta_1^2}\E_{(N, d_N)}\left(e^{2\tau_Y^N\wedge T}\right)&\leq \frac{1}{\eta_1^2}\max_{k\leq 10}\E_{(N, k)}\left(e^{2\tau_Y^N\wedge T}\right)\\
		&\leq \frac{\eps_0}{4} + \eta_1^{-2}\E_{0}\left(e^{2\tau_L\wedge T} \right)\\
		&\leq \frac{\eps_0}{4} + \frac{2+C_\tau\eta_1}{\eta_1^2},
	\end{align*}
	using first Proposition~\ref{PropScaling} and then Proposition~\ref{PropStopTime}.
	Combining the above bounds, we get for $N\geq N_0=\max\{N_{01}, N_{02}, N_{03}\}$ and $T>0$ that verifies Relation~\eqref{CondConst}, for $(N, d_N)\in \tilde{\cal{C}}_c$, 
	\[
		N^{-2}\E_{(N, d_N)}\left(\tilde{V}_C(R_Y^N(\tau_Y^N\wedge T), D_Y^N(\tau_Y^N\wedge T))\right)\leq \eps_0 +\frac{2+C_\tau\eta_1}{\eta_1^2}, 
	\]
	and since $N^{-2}\tilde{V}_C(N, d_N)\geq (1/\eta_0)^2$, by~\eqref{CondConst}, we can choose $\gamma>0$ small enough so that Relation~\eqref{EqPropStabY} holds, which concludes the proof of Proposition~\ref{PropStabY2}. 
\end{proof}

In the following Lemma, we show that the hitting time of $\cal{C}_c^c$ is monotomic in $d_N$, i.e., that the process takes more time to reach the boundary of $\cal{C}_c$ if it starts from a state on the diagonal $(D(0)=0)$ than if it starts away from it $(|D(0)|\geq 1)$. 

We present this result in a setting where $D(t)$ has transitions $\pm 1$, see Relation~\eqref{EqDLemma} to simplify the proof. The result can however easily be generalized to our exact process. 

\begin{lemma}\label{LemmaOrderTau}
	Let $N\geq 1$, and $d_I\in \N$ and $\alpha>0$. 
	Let $(R(t), D(t))$ be a continuous time Markov process starting from $(N, d_I)$,  and whose transitions are given by, for $r,d\in \N$, 
	\begin{equation}\label{EqDLemma}
		(r,d)\mapsto 
		\begin{cases}
			(r+1, d+1)\quad \text{with rate} \quad (r+\alpha d)/2\\
			(r+1, d-1)\quad \text{with rate} \quad \ind{d\geq 1} (r-\alpha d)/2 . 
		\end{cases}
	\end{equation}
	Introduce the stopping time 
	\[
		\tau\steq{def} \inf\{t\geq 0: D(t)\geq \eta_1\sqrt{R(t)}\}, 
	\]
	then for all $t\geq 0$, 
	\[
		\P_{(N, d_I)}\left(\tau\geq t \right)
		\leq \max\left\{\P_{(N, 0)}\left(\tau\geq t \right), \P_{(N, 1)}\left(\tau\geq t\right)\right\}.  
	\]
\end{lemma}

\begin{proof}
	Assume that $d_I=2d_0$ is even. The case of odd $d_I$ will be discussed at the end of the proof.
	We call $(R^0(t), D^0(t))$, resp. $(R^{d_0}(t), D^{d_0}(t))$ the process starting at $(N,0)$, resp. $(N, 2d_0)$. 
	Let $(t_i^0)_{i\geq 1}$, resp. $(t_i^{d_0})_{i\geq 1}$ be the increasing sequence of instant of jumps of the process $(R^0(t), D^0(t))$, resp. $(R^{d_0}(t), D^{d_0}(t))$: setting $t_0^a=0$ for $a\in \{0, d_0\}$, for all $i\geq 1$, 
	\[
		t_i^a=\inf\{t>t_{i-1}^a: \quad (R^a(t), D^a(t))\neq (R^a(t{-}), D^a(t{-}))\}. 
	\]
	Then, for all $i\in \N$, 
	\begin{align}\label{e:equalr}
		R^{d_0}(t_i^{d_0})=R^0(t_i^0)=N+i, 
	\end{align}
	and for $a\in \{d_0, 0\}$, the sequence $(E_i^a)\steq{def}(t_{i+1}^a-t_i^a)$ is constituted of independent random variables with $E_i^a$ distributed as an exponential random variable with mean $1/(N+i)$.
	We introduce the embedded Markov chain of the process $(D^{d_0}(t))$ and $(D^0(t))$: for all $i\in \N$, 
	\[
		A_i=D^{d_0}(t_i^{d_0}) \quad \text{and}\quad B_i=D^0(t_i^0).
	\]
	We will now show, using a coupling argument, that we can define random sequences $(\tilde{A}_i, \tilde{B}_i)$ with the same distribution as $(A_i, B_i)$ and for which almost surely, for all $i\in \N$, $\tilde{A}_i\geq \tilde{B}_i$.
	
	Note that $A_i$ and $B_i$ have the same parity as $i$. We chose their initial states with the same parity in order to prevent a ``crossing without meeting'' of the processes, that could happen around $0$: if $B_0=0$ and $A_0=1$, on the next step $B_{1}=1$ because of the reflexion at $0$, and $A_{1}$ could go down to $0$ with positive probability, which would contradict the inequality we want to show. This situation can't happen choosing both initial states $A_0$ and $B_0$ with the same parity. 
	
	We can construct both processes step by step, using some sequence $(U_i)$ of i.i.d. uniform variables on $[0,1]$. 
	For $i\geq 1$ and $x\in \N$, we introduce 
	\[
		\zeta_i(x)=
		\begin{cases}
			2\ind{U_i\leq \frac{N+i+\alpha x}{2(N+i)}}-1 \quad \text{if}\quad &x\geq 1\\
			1\quad \text{if} \quad &x=0. 
		\end{cases}
	\]
	Note that $\zeta_i(x)\in \{-1, 1\}$ for all $i\geq 0$, $x\in \N$. 
	Setting for all $i\geq 0$, 
	\[
		\tilde{A}_i= 2d_0+ \sum_{j=0}^{i-1}\zeta_j(\tilde{A}_j) \quad \text{and}\quad \tilde{B}_i= \sum_{j=0}^{i-1}\zeta_j(\tilde{B}_j), 
	\]
	it holds that $(\tilde{A}_i, i\geq 0)\steq{dist} (A_i, i\geq 0)$ and $(\tilde{B}_i, i\geq 0)\steq{dist} (B_i, i\geq 0)$. 
	
	Furthermore, for all $i\geq 0$, for $x\geq y$ such that $(x-y)\in 2\Z$, 
	\[
		x+\zeta_i(x)\geq y+\zeta_i(y) \quad \text{and} \quad (x+\zeta_i(x)-y-\zeta_i(z))\in 2\Z.  
	\]
	Therefore, since $\tilde{A}_0=2d_0\geq \tilde{B}_0=0$, for all $i\geq 1$, 
	\begin{equation}\label{IneqAB}
		\tilde{A}_i\geq \tilde{B}_i. 
	\end{equation}
	To conclude the proof, we now introduce the stopping times $V_{d_0}$ and $V_{0}$ as
	\begin{align*}
		V_{d_0}&\steq{def} \inf\{i\geq 0: \ \tilde{A_i}\geq \eta_1\sqrt{N+i}\}\\
		V_{0}&\steq{def} \inf\{i\geq 0: \  \tilde{B_i}\geq \eta_1\sqrt{N+i}\}. 
	\end{align*}
	Let $(E_i)$ be a sequence of independent random variables, where $E_i$ is exponentially distributed, with mean value $1/(N+i)$. Then, recalling Relation~\eqref{e:equalr} and that the distribution of $E_i^{\{0,d_0\}}$ only depends on the value of $R_i^{\{0,d_0\}}$, we have
	\[
		\tau_{|D(0)=2d_0}\steq{dist} \sum_{i=0}^{V_{d_0}} E_i \quad \text{and} \quad \tau_{|D(0)=0}\steq{dist} \sum_{i=0}^{V_0}E_i. 
	\]
	Using Relation~\eqref{IneqAB}, one has for all $k\geq 1$, 
	\[
		\{V_{d_0}\geq k\}\subset \{V_{0}\geq k\}, 
	\]
	and therefore, for $t\geq 0$, 
	\[
		\P(\tau\geq t\ |\ (R(0), D(0))=(N, 2d_0))\leq  \P(\tau\geq t\ |\ (R(0), D(0))=(N,0)), 
	\]
which concludes the proof for an even initial state. For an odd initial state, the same strategy can be used, comparing the process starting from $d_I=2d_0+1$ and the process starting from $1$. 	
\end{proof}

\subsection{Back to the initial process}\label{SecBack}

The results obtained in the previous sections have been shown on the simplified process $(Y(t))$, for which we only took into account the two main reactions. We now show that these results also apply to the process $(Z(t))$, whose $Q$-matrix is given by Relation~\eqref{QMatrixZ}. We assume that $(Y^N(t))$ and $(Z^N(t))$ are defined on the same probability space, and start form the same initial state $x_N\in \cal{C}_c$, with norm $\|x_N\|_1=N$ that can be chosen as large as necessary. 

The following lemma shows that until a time $t_1>0$, both processes stay at a distance of order $1$. It will only be used on trajectories stopped when they leave the interior cone $\cal{C}_c(\eta_1)$ defined in Relation~\eqref{DefIntCone}. To state our result rigorously, we introduce the following stopping time
\begin{equation}\label{EqTCN}
	T_C^N\steq{def} \inf\{t\geq 0: \ (Y^N(t),Z^N(t))\notin \cal{C}^2 \quad \text{or} \quad \min\{\|Z^N(t)\|_1, \|Y^N(t)\|_1\} \geq N/2\}
\end{equation}

\begin{lemma}\label{LemmaGronwall}
	Let $t_1>0$. We can find some $N_0\geq 1$ and $C_{t_1}>0$  such that for any initial state $x_N\in \cal{C}_c$ such that $\|x_N\|_1=N\geq N_0$, 
	$$\E_{x_N}\left(\sup_{t\leq t_1}\|Z(t\wedge T_C^N)-Y(t\wedge T_C^N)\|_1\right) <C_{t_1}.$$
\end{lemma}

\begin{proof}
We take $N_0$ large enough for the process $(Z^N(t))$ to be well defined until $T_C^N$. 
The proof relies on an argument using Gronwall's lemma. 
Using some Poisson processes $\cal{P}_i$ for $i=3,4,5$, and both $\cal{P}_1$ and $\cal{P}_2$ introduced in Relation~\eqref{SDEY}, the process $(Z(t))$ is solution of the SDE 
\begin{equation}\label{SDEZ}
	\begin{cases}
		\diff Z_1(t)= 2\cal{P}_1((0, Z_1(t{-})-2), \diff t)-\cal{P}_2((0, Z_2(t{-})-2), \diff t)\\
		\quad -4\cal{P}_3((0, \frac{(Z_1(t{-}))^{(4)}}{(Z_1(t{-}))^{(2)}(Y_2(t{-}))^{(2)}}), \diff t)+\cal{P}_5((0, \frac{1}{(Z_1(t{-}))^{(2)}(Y_2(t{-}))^{(2)}}), \diff t),\\
		\diff Z_2(t)= -\cal{P}_1((0, Z_1(t{-})-2), \diff t)+2\cal{P}_2((0, Z_2(t{-})-2), \diff t)\\
		\quad -4\cal{P}_4((0, \frac{(Z_2(t{-}))^{(4)}}{(Z_1(t{-}))^{(2)}(Y_2(t{-}))^{(2)}}), \diff t)+\cal{P}_5((0, \frac{1}{(Z_1(t{-}))^{(2)}(Y_2(t{-}))^{(2)}}), \diff t).
	\end{cases}
\end{equation}

We set for $i=1,2$, 
$$D_i(t)=|Z_i(t)-Y_i(t)| \quad \text{and}\quad V_i(t)=\min\{Z_i(t), Y_i(t)\}.$$

We only consider the processes $(Y(t))$ and $(Z(t))$ while they are in the cone $\cal{C}$. Therefore, for all $t\leq T_C^N$, 
$$pY_2(t)\leq Y_1(t), \quad pY_1(t)\leq Y_2(t), \quad pZ_2(t)\leq Z_1(t)\quad \text{and}\quad pZ_1(t)\leq Z_2(t),$$
where $p$ is given in Relation~\eqref{DefSets}. 
Integrating equations~\eqref{SDEZ} and~\eqref{SDEY}, for $t\leq T_C^N$, we have 
\begin{align*}
	D_1(t)&\leq \int_0^t2\cal{P}_1([V_1(s{-})-2, V_1(s{-})-2+D_1(s{-})), \diff s)\\
	&\quad +\int_0^t\cal{P}_2([V_2(s{-})-2, V_2(s{-})-2+D_2(s{-})), \diff s) \\
	&\quad + 4\int_0^t \cal{P}_3((0, \frac{(Z_1(t{-}))^{(4)}}{(Z_1(t{-}))^{(2)}(Y_2(t{-}))^{(2)}}), \diff t)+\cal{P}_5((0, \frac{1}{(Z_1(t{-}))^{(2)}(Y_2(t{-}))^{(2)}}), \diff t),\\
	&\leq \int_0^t2\cal{P}_1([V_1(s{-})-2, V_1(s{-})-2+\sup_{u\leq s}D_1(u{-})), \diff s)\\
	&\quad +\int_0^t\cal{P}_2([V_2(s{-})-2, V_2(s{-})-2+\sup_{u\leq s}D_2(s{-})), \diff s) \\
	&\quad + 4\int_0^t \cal{P}_3((0, \frac{2}{p^2}), \diff t)+\int_0^t\cal{P}_5((0, 1), \diff t),\\
\end{align*}
and therefore
\[
	\E\left(\sup_{s\leq t}D_1(s)\right)\leq \int_0^t 2\E\left(\sup_{u\leq s} D_1(u{-})\right)\diff s +\int_0^t \E\left(\sup_{u\leq s} D_2(u{-})\right)\diff s + \left(\frac{8}{p^2}+1\right)t. 
\]
A similar inequality can be shown for $D_2(t)$, which leads to 
\[
	\E\left(\sup_{s\leq t}\|Z(s)-Y(s)\|_1\right)\leq 6\int_0^t \E\left(\sup_{u\leq s} \|Z(u)-Y(u)\|_1\right)\diff s+ 2\left(\frac{8}{p^2}+1\right)t. 
\]
We conclude with Gronwall's Lemma, setting 
\[
	C_{t_1}= 2\left(\frac{8}{p^2}+1\right)t_1e^{6t_1}. 
\]
\end{proof}

\subsection{Proof of Proposition~\ref{p:filonov1}} We are now able to show the inequality~\eqref{InFilo} in the cone :

\begin{proof}
We first show the inequality for the process $(Z^N(t))$. 
We are going to show that we can find some $N_0\geq 1$ and some $\gamma>0$ such that for $x_N\in \cal{C}_c$ such that $\|x_N\|_1=N\geq N_0$,  
\begin{equation}\label{InFiloZ}
	N^{-2}\E_{x_N}\left(V_C(Z^N(\tau_Y^N\wedge T))-V_C(x_N)\right)\leq -\gamma, 
\end{equation}
where $\tau_Y^N$ has been defined in Relation~\eqref{DefTauY}. 

We define the process $(R_Z^N(t))=(Z_1^N(t)+Z_2^N(t))$.

It is straightforward that when $N$ is large, $T_C^N\geq \max\{\tau_Y^N\wedge T\}$ almost surely. In the following proof, since we look at the processes only until $\tau_Y^N$, we can forget about this stopping time when using Lemma~\ref{LemmaGronwall}. 

Using the same arguments as in the proof of Proposition~\ref{PropScaling}, we can show that for $T>0$, for $p\geq 1$, for the convergence in distribution, the limit 
\begin{equation}\label{SecGRZ}
	\lim_{N\to +\infty}\left(\frac{R_Z^N(t)^p}{N^p}\right)=(e^{pt})
\end{equation}
holds. 

Let $T>0$. We have for $N\geq 1$, and $x_N\in \cal{C}_c$ such that $\|x_N\|_1=N$, 
\begin{align*}
	N^{-2}&\E_{x_N}\left(V_C(Z^N(\tau_Y^N\wedge T))-V_C(x_N)\right)\\
	&= N^{-2}\E_{x_N}\left(V_C(Y^N(\tau_Y^N\wedge T))-V_C(x_N)\right)\\
	&+N^{-2}\E_{x_N}\left(V_C(Z^N(\tau_Y^N\wedge T))-V_C(Y^N(\tau_Y^N\wedge T))\right)\\
\end{align*}

Using Relation~\eqref{EqPropStabY} in Proposition~\ref{PropStabY2}, we can find some $N_{01}$ large enough so that $N^{-2}\E_{x_N}\left[V_C(Y^N(\tau_Y^N\wedge T))-V_C(x_N)\right]<-\gamma$ for some $\gamma>0$. This is the term useful here. We show that for $N$ large enough, the other term can be bounded by $\gamma/4$.  

\vspace{0.3cm}

We can find a constant $C_V>0$ such that for any $y,z\in \cal{C}$, 
\[
	|V(x)-V(y)|\leq C_V\|x-y\|_1(\|x\|_1^{3/2}+\|y\|_1^{3/2}).
\]

Therefore, 
\begin{align*}
	N^{-2}\E_{x_N}&\left(\left|V(Z^N(\tau_Y^N\wedge T))-V(Y^N(\tau_Y^N\wedge T))\right|\right)\\
	&\leq C_V \E_{x_N}\left(\frac{\sup_{s\leq T} \|Z^N(s)-Y^N(s)\|}{\sqrt{N}}\frac{R_Y^N(\tau_Y^N\wedge T)^{3/2}+R_Z^N(\tau_Y^N\wedge T)^{3/2}}{N^{3/2}}\right)\\
	&\leq C_V \E_{x_N}\left(\frac{(\sup_{s\leq T} \|Z^N(s)-Y^N(s)\|)^2}{N}\right)^{1/2}\\
	&\quad \quad \times\left(\E_{x_N}\left(\frac{R_Y^N(\tau_Y^N\wedge T)^3}{N^3}\right)^{1/2}+\E_{x_N}\left(\frac{R_Z^N(\tau_Y^N\wedge T)^3}{N^{3}}\right)^{1/2}\right)
\end{align*}
where the last inequality is due to Cauchy-Schwarz's inequality. The terms containing $R_Y$ and $R_Z$ are bounded because of the convergences mentioned in Proposition~\ref{PropScaling} and Relation~\eqref{SecGRZ}. Therefore, using Lemma~\ref{LemmaGronwall}, we can find some $N_{02}\geq 1$ such that 
\[
	N^{-2}\E_{x_N}\left(\left|V(Z^N(\tau_Y^N\wedge T))-V(Y^N(\tau_Y^N\wedge T))\right|\right)<\gamma/4, 
\]
which concludes the proof of Relation~\eqref{InFiloZ}. 

\vspace{0.8cm}

Now for the process $(X^N(t))$, we set the stopping time $\tau_X^N$ as 
\[
	\tau^N_X\steq{def} \inf\{t\geq 0: \lambda_t^N\geq (\tau_Y^N\wedge T)\}, 
\]
where $\lambda_t^N$ is defined in Relation~\eqref{EqLambda}. 

We have $X^N(\tau_X^N)=Z^N(\lambda_{\tau_X^N}^N)=Z^N(\tau_Y^N\wedge T)$. Besides, for $s>0$, 
\[
	\{\tau_X^N>s\}=\left\{\int_0^{\tau_Y^N\wedge T}\frac{1}{Z^N_1(u)^{(2)}Z^N_2(u)^{(2)}}\diff u>s\right\}.
\]
Until $\tau_Y^N\wedge T$, both $Z_1^N(t)$ and $Z^N_2(t)$ are larger than $N/2$, since $\tau_Y^N\wedge T\geq T_C^N$ with $T_C^N$ defined in Relation~\eqref{EqTCN}. Therefore, we have 
\begin{align*}
	\E(\tau_X^N)&=\int_{s=0}^{+\infty}\P\left(\tau_X^N>s\right)\diff s\\
		&\leq\int_{s=0}^{+\infty}\P\left(\int_0^{\tau_Y^N\wedge T}\frac{16}{(N-2)^4}\diff u>s\right)\diff s\\
		&\leq \int_{s=0}^{+\infty}\P\left(\tau_Y^N\wedge T> C N^4s\right)\\
		&\leq \frac{1}{C N^4}\E(\tau_Y^N\wedge T)
\end{align*}
where $C>0$ is a constant. Taking $N$ large enough, we get $\E_{x_N}(\tau_X^N)\leq 1$. All of this leads to Relation~\eqref{EqFiloGronwall}. 
\end{proof}
\section{Sides of the cone}\label{SecSides}

Remains to show that Relation~\eqref{InLyap} holds for $x\in \cal{C}_R\cup\cal{S}_R$, with $\|x\|_1$ large enough. 
Here, we use the Cartesian coordinates. The main difficulty here is located at the interface between the cone $\cal{C}$ and the sides $\cal{S}_R$ and $\cal{S}_L$, i.e., in the set $\Delta\steq{def}\{x\in \N^2: |x_2-px_1|\leq 5\}$, where the Lyapunov function changes of expression.

First, lets show the following :
\begin{lemma}\label{PropSide1}
	There exists $\eta_0>0$ in Relation~\eqref{DefIntCone} and $r^*>0$ such that
	\begin{itemize}
		\item[] for $x\in \cal{C}_R\cup\Delta$, $\|x\|\geq r^*$, the inequality $\cal{L}(V_C)(x)\leq -1$ holds. 
		\item[] for $x\in \cal{S}_R\cup\Delta$, $\|x\|\geq r^*$, the inequality $\cal{L}(V_R)(x)\leq -1$ holds. 
	\end{itemize}
\end{lemma}

\begin{proof}
	We start by the second inequality. Let $x\in \cal{S}_R\cup\Delta$. If $x_2\leq 2$, 
	\[
		\cal{L}(V_R)(x)=-4c_Sx_1^{(4)}\leq -1, 
	\]
	where $c_S$ is defined in Relation~\eqref{DefV}
	the last inequality holds if we take $r^*\geq (1/c_S)^{1/4}$. 
	Now we assume that $x_2\geq 2$. Since $2p x_1\geq x_2$ in $\cal{S}_R\cup \Delta$, $x_1$ can be chosen as large as necessary, so that 
	\begin{multline*}
		\cal{L}(V_R)(x)= c_Sx_1^{(2)}x_2^{(2)}(-4(x_1-2)+11(x_2-2))-4c_Sx_1^{(4)}-24c_Sx_2^{(4)} +7c_S,\\
				\leq -c_Sx_1^{(2)}x_2^{(2)}(4x_1-11x_2) +7c_S, 
	\end{multline*}
	and choosing $p<1/6$ and $r^*$ large enough allows us to conclude. 
	
	\vspace{0.5cm}
	
	To show the second inequality, we will use the axial coordinates.
	For $(r,d)\in \tilde{\cal{C}}_R$, one has
	\[
		\eta_0 \sqrt{r}<d<\frac{1-p}{1+p}r, 
	\]	
	by definition of $\tilde{\cal{C}}_c$ for the left inequality, see Relation~\eqref{DefIntCone}, and by definition of $\cal{S}_R$ for the right inequality, see Relation~\eqref{DefSets}. 

	To prove the lemma, we need the following inequalities. Using a Taylor expansion, we can find a function $g$ and some constants $C, C', r^*>0$ such that for $a, b\in [-5, 5]^2$, $(r, d)\in\tilde{\cal{C}_R}$,  
	\begin{equation}\label{Eq1Lemma9}
		\left|\frac{\tilde{V}_C(r+a, d+b)-\tilde{V}_C(r,d)}{r^3/d^2}-\left(\frac{3a}{r}-\frac{2b}{d}+\frac{3b^2}{d^2}\right)\right|\leq g(r,d), 
	\end{equation}
	and for $r\geq r^*$, 
	\begin{equation}\label{Eq2Lemma9}
		g(r,d)\leq \frac{C}{r^{2}}+\frac{C}{d^3},
	\end{equation}
	and, since the leading term in the round brackets in Relation~\ref{Eq1Lemma9} is $2b/d$, we have
	\begin{equation}\label{Eq3Lemma9}
		\left|\tilde{V}_C(r+a, d+b)-\tilde{V}_C(r,d)\right|\leq C'\frac{r^3}{d^3}. 
	\end{equation}
	
	We can distinguish the different contributions of the reactions to $\cal{A}(\tilde{V}_C)(r,d)$. We write for $(r,d)\in\tilde{\cal{C}_R}$ with $r\geq r^*$,
	\[
		\cal{A}(\tilde{V}_C)(r,d)=F_1(r,d) +F_2(r,d)+F_3(r,d)+F_4(r,d), 
	\]
	where the functions $F_i(r,d)$ are defined the following way:
	\begin{itemize}
		\item In the set $\tilde{\cal{C}}_R$, reactions $3A+2B\rightharpoonup 5A+B$ and $2A+3B\rightharpoonup A+5B$ dominate the evolution of the system. The contribution of the 3 other reactions is negligible in the value of $\cal{A}(\tilde{V}_C)(r,d)$. The term $F_4(r,d)$ gathers these contributions.
		\item Since $x_1=(r+d)/2$ and $x_2=(r-d)/2$ are both large, we can make the approximation $x_i^{(k)}\approx x_i^k$ for $i\in \{1,2\}$ and $k\geq 1$. $F_3(r,d)$ corresponds to the small correction due to this approximation. 
		\item The approximation~\eqref{Eq1Lemma9} also gives us a negligible term, gathered in $F_2(r,d)$.
		\item Finally, and most importantly, $F_1(r,d)$ contains the contribution of the two main reactions, after both approximations. 
	\end{itemize}
	
	To simplify the expressions, we will sometimes use the notation $(x_1, x_2)=\Phi^{-1}(r,d)=((r+d)/2, (r-d)/2)$. 
	
	The term $F_1(r,d)$ is given by:
	\begin{align*}
		F_1(r,d)&\steq{def} x_1^2x_2^2\frac{r^3}{d^2}\left(\frac{r+d}{2}\left(\frac{3}{r}-\frac{6}{d}+\frac{27}{d^2}\right)+\frac{r-d}{2}\left(\frac{3}{r}+\frac{6}{d}+\frac{27}{d^2}\right) \right)\\
		&=x_1^2x_2^2\frac{r^3}{d^2}\left(-3 +\frac{27r}{d^2}\right). 
	\end{align*}
	
	Using Relations~\eqref{Eq1Lemma9} and~\eqref{Eq2Lemma9}, the term $F_2(r,d)$ can be upper bounded as follow, 
	\[
		|F_2(r,d)|\leq x_1^2x_2^2\frac{r^3}{d^2}\left(\frac{C}{r^2}+\frac{C}{d^3}\right)\\
							\leq \frac{r^7}{d^4}, 
	\]
	where the last inequality holds if $r^*$ is chosen large enough. 
	
	The term $F_3(r,d)$ is defined as
	\begin{multline*}
		F_3(r,d)\steq{def} (x_1^{(3)}x_2^{(2)}-x_1^3x_2^2)(\tilde{V}_C(r+1, d+3)-\tilde{V}_C(r, d)) \\
		+(x_1^{(2)}x_2^{(3)}-x_1^2x_2^3)(\tilde{V}_C(r+1, d-3)-\tilde{V}_C(r, d)), 
	\end{multline*}
	which can be upper bounded, using Relation~\eqref{Eq3Lemma9}, as follows
	\[
		|F_3(r,d)|\leq C_3 r^4 C'\frac{r^3}{d^3}= C_3'\frac{r^7}{d^3}
	\]
	for some $C_3, C_3'>0$. 
	
	Finally, $F_4(r,d)$ is defined as follows: 
	\begin{multline*}
		F_4(r,d)\steq{def} x_1^{(4)}(\tilde{V}_C(r-4, d-4)-\tilde{V}_C(r, d)) +x_2^{(4)}(\tilde{V}_C(r-4, d+4)-\tilde{V}_C(r, d))\\
		+(\tilde{V}_C(r+2, d)-\tilde{V}_C(r, d)), 
	\end{multline*}
	and can be upper bounded, using Relation~\eqref{Eq3Lemma9}, as follows
	\[
		\left|F_4(r,d)\right|\leq C_4 \frac{r^7}{d^3}. 
	\]
	
	Combining the above bounds, we get the following inequality, for $(r,d)\in\tilde{\cal{C}_R}$ with $r\geq r^*$,
	\[
		\cal{A}(\tilde{V}_C)(r,d)\leq x_1^2x_2^2\frac{r^3}{d^2}\left(-3 +\frac{27r}{d^2}\right) +(C_3'+C_4)\frac{r^7}{d^3}+\frac{r^7}{d^4}, 
	\]
	and choosing $\eta_0>3$, and $r^*$ large enough, it is easy to conclude that $\cal{A}(\tilde{V}_C)(r,d)=\cal{L}(V_C)(x_1, x_2)\leq -1$, when $(r,d)\in \tilde{\cal{C}}_R$, and $(x_1, x_2)=((r+d)/2, (r-d)/2)\in \cal{C}_R$. Since the choice of $p$ is not restrictive here, it is straightforward that this inequality holds also for the states in $\Delta$.  
\end{proof}

We now establish our bound at the interface between $\cal{C}$ and $\cal{S}$, i.e., where the definition of $V$ changes. If the process is at the state $x\notin \Delta$, all reactions keep the process on the same side of the border $\{x_2=px_1\}$. Therefore, for $x\notin \Delta$, 
\[
	\begin{cases}
		\cal{L}(V)(x)=\cal{L}(V_C)(x) \quad \text{if}\quad x\in \cal{C}_R, \\
		\cal{L}(V)(x)=\cal{L}(V_R)(x) \quad \text{if}\quad x\in \cal{S}_R.
	\end{cases}
\]

However, for $x\in \Delta$, one has to proceed more carefully. 

\begin{lemma}\label{PropSide2}
	If $p<1/29$, for $x\in \Delta$ such that $\|x\|\geq r^*$, for $r^*$ large enough, 
	\begin{equation}\label{InBorderDelta}
		\begin{cases}
			V_C(x)\leq V_R(x) \quad \text{if}\quad x\in \cal{C}_R, \\
			V_C(x)\geq V_R(x) \quad \text{if}\quad x\in \cal{S}_R.
		\end{cases}
	\end{equation}
	As a consequence, the inequality $\cal{L}(V)(x)\leq-1$ holds.
\end{lemma}
\begin{proof}
	First, we assume that Relation~\eqref{InBorderDelta} holds. For $x\in \Delta\cap \cal{C}$, if after a reaction the process reaches the state $y\in \cal{S}_R$, we have: 
	\[
		V(y)-V(x)=V_R(y)-V_C(x)\leq V_C(y)-V_C(x), 
	\]
	and therefore
	\[
		 \cal{L}(V)(x)\leq \cal{L}(V_C)(x)\leq -1, 
	\]
	the last inequality being a consequence of Lemma~\ref{PropSide1}. A similar argument can be used for $x\in \Delta\cap \cal{S}_R$. 
	
	Now we prove both inequalities. Let $x_1\geq r^*/2$, and $a\in [-5,5]^2$, $\|a\|\geq 1$. We set $x=(x_1, px_1)$. 
	We can find some constant $C>0$ such that  
	\[
		\left|V_C(x+a)-V_C(x)-\croc{\nabla V_C(x),a}\right|\leq \frac{C}{x_1},
	\]
	where $\nabla V_C(x)$ is the gradient of $V_C$ at $x$: 
	\[
		\nabla V_C(x)=\frac{(1+p)^2}{(1-p)^3}\binom{1-5p}{5-p}. 
	\]
	Furthermore, we have that  
	\[
		V_R(x+a)-V_R(x)=c_S(a_1+6a_2)=\croc{\nabla V_R(x),a}
	\]
	where $\nabla V_R(x)$ is the gradient of $V_R$ at $x$: 
	\[
		\nabla V_R(x)=c_S\binom{1}{6}\quad  \text{with} \quad c_S=\frac{(1+p)^3}{(1-p)^2(1+6p)}.
	\]
	Since by definition of $c_S$, $V_R(x)=V_C(x)$, we get 
	\[
		\left|V_C(x+a)-V_R(x+a)-(\nabla V_C(x)-\nabla V_R(x))\cdot a\right|\leq \frac{C}{x_1}, 
	\]
	and after calculation, we can show that  
	\[
		(\nabla V_C(x)-\nabla V_R(x))\cdot a=\frac{(1+p)^2(1-29p)}{(1-p)^3(6p+1)} (pa_1-a_2). 
	\]
	Choosing $p<1/29$, and $r^*$ large enough, we conclude that $V_R(x+a)-V_C(x+a)$ has the same sign as $(pa_1-a_2)$, and since $x+a\in \cal{C}$ if $pa_1<a_2$, and $x+a\in \cal{S}_R$ if $pa_1>a_2$, this proves the claim. 
\end{proof}

\section*{Acknowledgements}
The authors thank Jonathan C. Mattingly and Philippe Robert for helpful discussions. AA acknowledges partial support from PRIN 2022 project ConStRAINeD (2022XRWY7W).

\bibliographystyle{plain}
\bibliography{ref}

\end{document}